\documentstyle[amstex,amssymb,11pt]{article}
\renewcommand{\ni}{\noindent}

\newcommand{\integ}[2]{\displaystyle \int_{#1}^{#2}}

\newtheorem{theorem}{Theorem}[section]
\newtheorem{proposition}{Proposition}[section]

\newtheorem{remark}{Remark}[section]

\def \R{I\!\!R}

\def \E{I\!\!E}

\def \bf{\textbf}
\def \it{\textit}

\def \no {\noindent}

\begin{document}
\author{Said
Hamad\`ene\thanks{Universit\'e du Maine, LMM, Equipe Statistique et
Processus, Avenue Olivier Messiaen, 72085 Le Mans, Cedex 9, France.
e-mail: hamadeneATuniv-lemans.fr}\,\, and \,\,Youssef Ouknine
\thanks{D\'{e}partement de Math\'{e}matiques,
Facult\'{e} des Sciences Semlalia, Universit\'{e} Cadi Ayyad,
Marrakech, Maroc. e-mail: ouknineATucam.ac.ma}}
\title{Reflected Backward SDEs with General Jumps}
\date{\today}
\maketitle
\begin{abstract}
In the first part of this paper we give a solution for the
one-dimensional reflected backward stochastic differential equation
(BSDE for short) when the noise is driven by a Brownian motion and
an independent Poisson random measure. The reflecting process is
right continuous with left limits (RCLL for short) whose jumps are
arbitrary. We first prove existence and uniqueness of the solution
for a specific coefficient in using a method based on a combination
of penalization and the Snell envelope theory. To show the general
result we use a fixed point argument in an appropriate space. The
second part of the paper is related to BSDEs with two reflecting
barriers. Once more we prove existence and uniqueness of the
solution of the BSDE.
\end{abstract}
\bigskip

\ni \bf{ Key Words:} Backward SDEs ; Penalization ; Poisson point
process ; Snell envelope ; Mokobodski's hypothesis.
\medskip

\ni \bf{AMS Classification} (1991): 60H10, 60H20, 60H99.
\section{Introduction}
Non-linear backward stochastic differential equations (BSDEs in
short) were introduced by Pardoux $\&$ Peng \cite{pp1} when the
noise is driven by a Brownian motion. The objective is to give a
probabilistic interpretation of a solution of a second order
quasi-linear partial differential equation. Since then, these
equations have gradually became an important mathematical tool which
is encountered in many fields of mathematics such as finance,
stochastic optimal control and games, partial differential equations
and so on (see e.g. \cite{bbp, CK, ekp, [EPQ], [H], [Pa], [PP1], tl}
and the references therein).

Later Tang $\&$ Li \cite{tl}, considered standard BSDEs when the
noise is driven not only by a Brownian motion but also by an
independent Poisson random measure. They showed existence and
uniqueness of the solution. Barles et al. \cite{bbp} studied the
link of those BSDEs with viscosity solutions of integral-partial
differential equations.

One barrier reflected BSDEs have been introduced by El-Karoui et al.
in \cite{ekp}. In their setting, one of the components of the
solution is forced to stay above a given barrier which is a
continuous adapted stochastic process. The main motivation in
\cite{ekp} is the pricing of American options especially in
constrained markets (see also \cite{[EPQ]}). The generalization to
the case of two reflecting barriers has been carried out by Cvitanic
$\&$ Karatzas in \cite{CK}.

Later, on the one hand, Hamad\`ene $\&$ Oukine (\cite{ho}) have
studied one reflecting barrier BSDEs when the noise is driven by a
Brownian motion and an independent Poisson measure. They showed
existence and uniqueness of the solution when the reflecting barrier
has only inaccessible jumps, i.e., jumps which stem only from the
Poisson part. On the other hand, S.Hamad\`ene \cite{[H]} has
introduced BSDEs with one right continuous with left limits
reflecting barrier in the case of Brownian noise. Since then there
have been several works on BSDEs with discontinuous barriers when
the noise comes only from a Brownian motion (\cite{lepeltiermingyu,
pengmingyu}).
\medskip

Therefore the main objective of this paper is to deal with reflected
BSDEs when the noise comes from a Brownian motion and an independent
Poisson process and the reflecting processes are just RCLL. No more
conditions are imposed on their jumps as it was e.g. in \cite{ho}.
They can be predictable or inaccessible. In our study we consider
the case of one reflecting barrier as well as the case of two
reflecting barriers. For both cases we show existence and uniqueness
of the solution when the coefficients of the BSDEs are Lipschitz.
\medskip

This work completes the known results on the same subject since the
jumps of the reflecting processes are arbitrary and the sources of
noise are twice, Brownian and Poisson. A second motivation of our
work is that the two barrier reflected BSDEs we consider here are
much involved in finance, especially when we deal with convertible
bonds in defautable markets (see e.g. \cite{jeanbl}). Finally, this
work opens a window towards viscosity solutions of variational
inequalities with discontinuous obstacles.
\medskip

This paper is organized as follows.
\medskip

Section 2 contains hypotheses and the setting of the problem. In
Section 3, we show uniqueness of the solution of the BSDE with one
reflecting RCLL barrier $(S_t)_{t\leq 1}$, a Lipschitz coefficient
$f(t,\omega,y,z,v)$ and square integrable terminal value $\xi$. In
Section 4, we address the question of existence of the solution for
the BSDE. Since there is a lack of comparison of solutions of
standard BSDEs whose noise contains a Poisson part, especially when
the coefficients depend on $v$ (see \cite{bbp} for a
counterexample), we first assume that $f$ does not depend on this
latter variable. Then we show the existence and uniqueness of the
solution of the BSDE. The method we used combines penalization with
the general theory of Snell envelope of processes. Later, in order
to obtain the result for general coefficients $f(t,\omega,y,z,v)$ we
have introduced a contraction in an appropriate Banach space of
processes which then has a fixed point which provides the unique
solution of the BSDE. At the end of this section we focus on some
properties of the reflecting process $K$ of the solution. Finally in
the last section, we address the problem with two reflecting
barriers. Once more, under the well-known Mokobodski's hypothesis,
we show existence and uniqueness of the solution.
\section{Setting of the problem and hypotheses} First a simplicity reason, we
fix the horizon $T$ of the problem  equal to $1$ and of course our
results still valid if $T\neq 1$. Let $(\Omega, {\cal F}, ({\cal
F}_t)_{t\leq 1})$ be a stochastic basis such that ${\cal F}_0$
contains all $P$-null sets of ${\cal F}$ and ${\cal
F}_{t+}=\bigcap_{\epsilon >0}{\cal F}_{t+\epsilon}={\cal F}_t$,
$\forall t< 1$. We assume that the filtration is generated by the
two following mutually independent processes:

- a $d$-dimensional Brownian motion $(B_t)_{t\leq 1}$,

- a random Poisson measure $\mu $ on $I\!\!R^{+}\times U$, where $%
U:=I\!\!R^{l}\setminus \{0\}$ is equipped with its Borel
$\sigma$-algebra ${\cal U}$,
with compensator $\nu (dt,de)=dt\lambda (de)$, such that $\{\tilde{\mu}%
([0,t]\times A)=(\mu -\nu )([0,t]\times A)\}_{t\leq 1}$ is a
martingale for every $A\in {\cal U}$ satisfying $\lambda (A)<\infty
$. The measure $\lambda $ is assumed to be $\sigma $-finite on
$(U,{\cal U)}$ and integrates the function $e\in U\mapsto 1\wedge
|e|^2$.
\medskip

Let us now introduce the followings:
\medskip

- ${\mbox{ L}}^{2}(U,d\lambda)$ is the space of deterministic Borel
functions $(\varpi(e))_{e\in U}$ from $U$ to $R$ which are
$d\lambda$-square integrable;

- ${\cal S}^2$ the set of ${\cal F}_t$-adapted RCLL processes
$(Y_t)_{t\leq 1}$ with values in $I\!\!R$ and $I\!\!E[\sup_{t\leq
1}|Y_t|^2]<\infty$. We denote by ${\cal S}^2_i$ the subset of ${\cal
S}^2$ which contains non-decreasing processes $K:=(K_t)_{t\leq 1}$
with $K_0=0$;

- ${\cal P}$ the $\sigma$-algebra of ${\cal F}_t$-progressively
measurable sets on $\Omega\times [0,1]$ and ${\cal H}^{2,k}$ the set
of ${\cal P}$-measurable processes $Z:=(Z_t)_{t\leq 1}$ with values
in $I\!\!R^k$ and $dP\otimes dt$-square integrable;

- ${\cal P}^d$ the $\sigma$-algebra of ${\cal F}_t$-predictable sets
on $\Omega \times [0,1]$ and ${\cal L}^{2}$ the set of mappings
$V:\Omega \times [0,1]\times U\rightarrow I\!\!R$, ${\cal
P}^d\otimes {\cal U}$-measurable and $dP\otimes dt\otimes
d\lambda$-square integrable;

- ${\cal T}$ the set of ${\cal F}_t$-stopping times with values in
$[0,1]$;

- the class [D] is the set of ${\cal F}_t$-adapted RCLL processes
$\zeta=(\zeta_t)_{t\leq 1}$ such that the set of random variables
$\{\zeta_\tau, \,\,\tau \in \cal T \}$ is uniformly integrable;

- for a given RCLL process $(w_t)_{t\leq 1}$, for any $t\leq 1$,
$w_{t-}=\lim_{s\nearrow t}w_s$ $(w_{0-}=w_0)$,
$\Delta_tw=w_t-w_{t-}$ and $w_-:=(w_{t-})_{t\leq 1}$.
\bigskip

We are now given three objects:
\medskip

\ni - a terminal value $\xi \in L^{2}(\Omega ,F_{1},P)$
\medskip

\ni - a map $f:\Omega \times \left[ 0,1\right] \times
I\!\!R^{1+d}\times {\mbox{ L}}^{2}(U,d\lambda)\longrightarrow
I\!\!R$
which with $(t,\omega ,y,z,v)$ associates $f(t,\omega ,y,z,v)$ and which is $%
{\cal P}\otimes {\cal B}(I\!\!R^{1+d})\otimes {\cal B}({\mbox{
L}}^{2}(U,d\lambda))$-measurable. In addition we assume that:
\medskip

$(i)$ the process $(f(t,0,0,0))_{t\leq 1}$ belongs to ${\cal
H}^{2,1}$ \medskip

$(ii)$ $f$ is uniformly Lipschitz with respect to $(y,z,v)$, i.e.,
there exists a constant $C_f\geq 0$ such that for any $(y,z)$,
$(y',z')$ $\in \R\times \R^{d}$ and $v,\,v'\in {\cal
L}^{2}_\lambda(U)$ we have:
\[
P-a.s.,\,\,|f(\omega ,t,y,z,v)-f(\omega ,t,y',z',v')|\leq
C_f(|y-y'|+|z-z'|+\|v-v'\|).
\]
- an "obstacle" process $S:=(S_t)_{t\leq 1}$ which belongs to ${\cal
S}^2$.
\bigskip

Let us now introduce the reflected BSDE with general jumps
associated with $(f,\xi ,S)$. A solution is a quadruple
$(Y,Z,K,V):=(Y_{t},Z_{t},K_{t},V_{t})_{t\leq 1}$ of processes with
values in $I\!\!R^{1+d}\times I\!\!R^{+}\times {\mbox{
L}}^{2}(U,d\lambda)$ such that:
\begin{equation}\label{prince}
\left\{
\begin{array}{ll}
(i)& Y\in {\cal S}^{2},\,Z\in {\cal H}^{2,d},\,V\in {\cal
L}^{2}\mbox{ and
}K\in {\cal S}^{2}_i \\
(ii)& Y_{t}=\xi +\displaystyle
\int_{t}^{1}f(s,Y_{s},Z_{s},V_{s})ds+K_{1}-K_{t}-
\int_{t}^{1}Z_{s}dB_{s}-\displaystyle \int_{t}^{1}\!\!\!\!
\int_{U}^{{}}V_{s}(e)\tilde{\mu}(ds,de)\,,\forall t\leq 1 \\
(iii)& Y\geq
S\\
(iv)& \mbox{if $K^c$ (resp. $K^d$) is the continuous (resp. purely
discontinuous) part of $K$, then}\\& \mbox{$K^d$ is ${\cal
P}^d$-measurable, }
\int_{0}^{1}(Y_{t}-S_{t})dK_{t}^{c}=0\,\,\text{and }\forall t\leq 1,
\Delta K_t^d=(S_{t-}-Y_t)^{+}1_{[ Y_{t-}=S_{t-}]}.
\end{array}
\right.
\end{equation}

The main reason for the second part of $(iv)$ is that the process
$Y$ has two types of jumps. The inaccessible ones which stem from
the Poisson martingale part $(\int_{0}^{t}\!\!
\int_{U}^{{}}V_s(e)\tilde{\mu}(ds,de))_{t\leq 1}$ and the
predictable ones which come from the predictable negative jumps of
$S$. Those latter are the source of the predictable jumps of $Y$ and
then also of $K$, which of course are the same. Thus the condition
$\Delta K_t^d=(S_{t-}-Y_t)^{+}1_{[ Y_{t-}=S_{t-}]}$ is just a
characterization of the predictable jumps of $Y$.

\begin{remark}: \label{rem1} The second part of condition
$(iv)$ implies in particular that
\begin{equation}
\label{eqrem} \mbox{$\int_{0}^{1}(Y_{t-}-S_{t-})dK_{t}=0.$}
\end{equation} Actually
$$
\begin{array}{ll}
\int_{0}^{1}(Y_{t-}-S_{t-})dK_{t}&=\int_{0}^{1}(Y_{t-}-S_{t-})dK_{t}^c+\int_{0}^{1}(Y_{t-}-S_{t-})dK_{t}^d\\
{}&=\int_{0}^{1}(Y_{t}-S_{t})dK_{t}^c+\sum_{t\leq
1}(Y_{t-}-S_{t-})\Delta K_{t}^d \\{}&=0.
\end{array}
$$
The last term of the second equality is null since $K^d$ jumps only
when $Y_{t-}=S_{t-}$. $\Box$
\end{remark}

To begin with, we are going to focus on the uniqueness of the
solution of the BSDE (\ref{prince}).
\section{Uniqueness of the solution}
\begin{proposition}\label{unicite}: Under the above assumptions on $f$, $\xi $
and $(S_{t})_{t\leq 1}$, the reflected BSDE (\ref{prince})
associated with $(f,\xi ,S)$ has at most one solution.
\end{proposition}

\ni $Proof$: Let us consider two solutions $(Y,Z,K,V)$ and
$(Y^{\prime },Z^{\prime },K^{\prime },V^{\prime })$ of
$(\ref{prince})$. First let us assume that for any $t\leq 1$ we have
$\left( Y_{t-}-Y_{t-}^{\prime }\right) \left( dK_{t}-dK_{t}^{\prime
}\right) \leq 0.$ Next for $t\leq 1$ let us set:
$$\begin{array}{c}
\Delta _{t}:=| Y_{t}-Y_{t}^{\prime }|
^{2}+\int_{t}^{1}|Z_{s}-Z_{s}^{\prime
}|^{2}ds+\int_{t}^{1}\!\!\!\int_{U}( V_{s}( e) -V_{s}^{\prime }( e)
) ^{2}\lambda (de)ds. \end{array}$$Then applying It\^{o}'s formula
with $(Y-Y^{\prime })^2$ and taking expectation yield:
$$\begin{array}{lll} I\!\!E[\Delta _{t}] & \leq &
2\int_{t}^{1}I\!\!E[\left( Y_{s}-Y_{s}^{\prime }\right) \left(
f(s,Y_{s},Z_{s},V_{s})-f(s,Y_{s}^{\prime },Z_{s}^{\prime
},V_{s}^{\prime })\right) ]ds, \,\,\forall t\leq 1.
\end{array}
$$
Therefore because of the inequality $2|ab|\leq 2a^2+\frac{b^2}{2}$
for any $a,b\in \R$ and since $f$ is uniformly Lipschitz we obtain:
$\forall \,\,t\leq 1$,
\begin{equation}
\label{eqdelta}
\begin{array}{ll}
{}&\E[ \left| Y_{t}-Y_{t}^{\prime }\right| ^{2}+\frac{1}{2}%
\int_{t}^{1}|Z_{s}-Z_{s}^{\prime }|^{2}ds+\frac{1}{2}
\int_{t}^{1}\int_{U}\left( V_{s}\left( e\right) -V_{s}^{\prime
}\left( e\right) \right) ^{2}\lambda \left( de\right) ds]\\
{}&\,\,\,\qquad \qquad\qquad \leq (2C_f+4C_f^2)\E[
\int_{t}^{1}\left( Y_{s}-Y_{s}^{\prime }\right) ^{2}ds].
\end{array}
\end{equation}
Then from Gronwall's lemma and the right continuity of $Y-Y^{\prime
}$ we get $Y=Y^{\prime }$. It follows from (\ref{eqdelta}) that
$\left( Z,V\right) =\left(Z^{\prime },V^{\prime }\right)$ and
finally $K=K'$. Whence the uniqueness of the solution
of (%
\ref{prince}).
\medskip

To complete the proof it remains to show that

\begin{equation}\label{inegalite}\int_{]t,1]}(Y_{s-}-Y'_{s-})(dK_{s}-dK_{s}^{\prime})\leq 0,
\forall t\leq 1.\end{equation} Actually for any $t\leq 1$ we have,
\begin{equation}\label{eqine} \int_{]t,1]}^{{}}(Y_{s-}-Y'_{s-})(dK_{s}-dK_{s}^{\prime})=
\int_{]t,1]}^{{}}(Y_{s-}-Y'_{s-})(dK_{s}^{c}-dK_s^{\prime c})+
\int_{]t,1]}^{{}}(Y_{s-}-Y'_{s-})(dK_{s}^{d}-dK_s^{\prime d}).
\end{equation} But as $Y$ and $Y'$ belong to ${\cal S}^{2}$ and
their jumps $\delta (\omega ):=\{t\in [0,1],\Delta _{t}Y(\omega)\neq
0\}$ and $\delta ^{\prime }(\omega ):=\{t\in [0,1],\Delta
_{t}Y^{\prime }(\omega)\neq 0\}$ are at most countable then:
\begin{equation}
\label{eqineg}
\begin{array}{ll}
&  \int_{]t,1]}^{{}}(Y_{s-}-Y'_{s-})(dK_{s}^{c}-dK_s^{\prime c})=
\int_{]t,1]}^{{}}(Y_{s}-Y_{s}^{\prime })(dK_{s}^{c}-dK_s^{\prime
c})=
\int_{]t,1]}^{{}}(Y_{s}-S_{s})(dK_{s}^{c}-dK_s^{\prime c}) \\
& \qquad + \int_{]t,1]}^{{}}(S_{s}-Y_{s}^{\prime
})(dK_{s}^{c}-dK_s^{\prime c})=-
\int_{]t,1]}^{{}}(Y_{s}-S_{s})dK_s^{\prime c}+
\int_{]t,1]}^{{}}(S_{s}-Y_{s}^{\prime })dK_{s}^{c}\leq 0.
\end{array}
\end{equation}
Additionally we have:
\begin{equation}
 \int_{]t,1]}^{{}}(Y_{s-}-Y'_{s-})(dK_{s}^{d}-dK_s^{\prime
d})= \int_{]t,1]}^{{}}(Y_{s-}-Y'_{s-})dK_{s}^{d}-
\int_{]t,1]}^{{}}(Y_{s-}-Y'_{s-})dK_s^{\prime d}.
\end{equation}
However
\begin{equation}
\label{eqk2}  \int_{]t,1]}^{{}}(Y_{s-}-Y'_{s-})dK_s^{\prime d}=
\int_{]t,1]}(Y_{s-}-S_{s-})dK_s^{\prime d}\geq 0
\end{equation}
since the jumps of ${K'}^d$ occur only when $Y'_-=S_-$. In the same
way we have $\int_{]t,1]}^{{}}(Y_{s-}-Y'_{s-})dK_s^{ d}\leq 0$ and
then $\int_{]t,1]}^{{}}(Y_{s-}-Y'_{s-})(dK_{s}^{d}-dK_s^{\prime
d})\leq 0$. This inequality and (\ref{eqineg}) lead to
(\ref{inegalite}). Thus the proof of uniqueness of the solution of
the BSDE (\ref{prince}) is complete. $\Box$
\bigskip

\section{Existence of the solution}
We are now going to prove that equation (\ref{prince}) has a
solution.  Our method combines penalization and the Snell envelope
method. Penalization as it has been used e.g. in \cite{ekp} cannot
be applied because of lack of comparison of the solutions of the
penalization scheme. To get this comparison we are led to assume, in
a first step, that the function $f(t,\omega,y,z,v)$ does not depend
on $v$ and for the sake of simplicity we will suppose moreover that
$f$ does not depend also on $(y,z)$, $i.e.$,
$f(t,\omega,y,z,v)=g(t,\omega)$. Later to obtain the result in the
general setting we will use a fixed point argument with an
appropriate mapping. Finally note that as a by-product of the
penalization method is the approximation of the solution of the
reflected equation by solutions of standard BSDEs, i.e. without
reflection. This point could be important when we deal with the
issue of either numerical schemes or the solution of BSDE
(\ref{prince}) or viscosity solutions of the related PDIEs.
\medskip

To begin with let us assume that the function $f$ does not depend on
$(y,z,v)$, i.e., P-$a.s.$,\\ $f(t,\omega ,y,z,v)\equiv g(t,\omega
)$, for any $t,y,z$ and $v$. In the following result, we establish
the existence of the solution of the BSDE associated with $(g,\xi
,S)$ in using the penalization method.
\begin{theorem}: The one barrier reflected BSDE (\ref{prince})
associated with $(g,\xi,S)$ has a unique solution
$(Y_t,Z_t,K_t,V_t)_{t\leq 1}$.
\end{theorem}
$Proof$: For $n\geq 0$, let $\left(
Y_{t}^{n},\;Z_{t}^{n},V_{t}^{n}\right)_{t\leq 1}$ be the ${\cal
F}_{t}$-adapted process with values in $I\!\!R^{1+d}\times
\mbox{L}^2(U,d\lambda)$, unique solution of the BSDE associated with
$(g(t,\omega)+n(y-S_t)^-,\xi)$ ($x^-:=max(0,-x)$, $\forall x\in R$),
which exists according to the results either by Tang $\&$ Li
\cite{tl} or Barles et al. \cite{bbp}, $i.e.$,
\begin{equation}  \label{esti1}
\left\{
\begin{array}{l}
Y^n\in {\cal S}^2, Z^n\in {\cal H}^{2,d}\mbox{ and }V^n\in {\cal
L}^2\\
Y_{t}^{n}=\xi +\int_{t}^{1}g(s)ds + \int_{t}^{1}n\left(
Y_{s}^{n}-S_{s}\right) ^{-}ds-\int_{t}^{1}Z_{s}^{n}dB_{s}
-\int_{t}^{1}\!\!\int_{U}V_{s}^{n}\left( e\right) \tilde{\mu }\left(
ds,de\right), \forall t\leq 1. \end{array} \right.\end{equation}

From now on the proof will be divided into four steps.
\medskip

\ni \bf{ Step 1}: For any $n\geq 0$, $Y^n\leq Y^{n+1}$. \medskip

For any $t\leq 1$, we have:
$$
\begin{array}{l}
Y^n_t-Y^{n+1}_t=\int_t^1(Z^n_s-Z^{n+1}_s)dB_s+\int_t^1\int_U(V^n_s(e)-V^{n+1}_s(e))\tilde{\mu}(ds,de)
\\\qquad\qquad\qquad\qquad+\int_t^1\{n(Y^n_s-S_s)^--(n+1)(Y^{n+1}_s-S_s)^-\}ds.
\end{array}
$$
But
$n(Y^n_s-S_s)^--(n+1)(Y^{n+1}_s-S_s)^-=b^n_s+a^n_s(Y^n_s-Y^{n+1}_s)$
where $b^n_s\leq 0$ and $|a^n_s|\leq n+1$, $\forall s\leq 1$. So let
us set $\Theta ^n_t=e^{\int_0^ta^n_sds}$, $t\leq 1$ ; then using
It\^o's formula we obtain:
$$
d(\Theta^n_t (Y^n_t-Y^{n+1}_t))=\Theta_t^n b^n_tdt+dM^n_t, \,\,t\leq
1$$where $(M_t^n)_{t\leq T}$ is a martingale. Taking now the
conditional expectation w.r.t. ${\cal F}_t$ we obtain $Y^n\leq
Y^{n+1}$ since $Y^{n+1}_1-Y^n_1=0$, $\Theta^n\geq 0$ and $b^n\leq
0$. $\Box$
\medskip

\ni \bf{ Step 2}: For any $n\geq 0$, the process $Y^n$ satisfies:
$$
\forall t\leq 1,\,\, Y_t^n=\mbox{ esssup}_{\tau \geq
t}\E[\int_t^\tau g(s)ds +(Y^n_\tau \wedge S_\tau) 1_{[\tau<1]}+\xi
1_{[\tau=1]}|{\cal F}_t].$$ Actually for any $n \geq 0$ and $t\leq
1$ we have:
\begin{equation}
\label{eqstep4} Y_{t}^{n}=\xi +\int_{t}^{1}g(s)ds
-\int_{t}^{1}Z_{s}^{n}dB_{s} + \int_{t}^{1}n\left(
Y_{s}^{n}-S_{s}\right) ^{-}ds
-\int_{t}^{1}\!\!\int_{U}V_{s}^{n}\left( e\right) \tilde{\mu }\left(
ds,de\right).\end{equation} Therefore for any stopping time $\tau
\geq t$ we have:
\begin{equation}
\label{eqyn1}
\begin{array}{ll} Y_{t}^{n}&=\E[Y_{\tau}^{n}+ \int_{t}^{\tau}g(s)ds
+\int_{t}^{\tau}n\left(
Y_{s}^{n}-S_{s}\right) ^{-}ds|{\cal F}_t]\\
{}&\geq \E[(S_\tau \wedge Y^n_\tau) 1_{[\tau<1]}+\xi 1_{[\tau=1]}+
\int_{t}^{\tau}g(s)ds |{\cal F}_t]\\
\end{array}\end{equation}
since $Y_{\tau}^{n}\geq (S_\tau \wedge Y^n_\tau) 1_{[\tau<1]}+\xi
1_{[\tau=1]}$. On the other hand let $\tau^*_t$ be the stopping time
defined as follows:
$$
\tau^*_t=\inf\{s\geq t,K^n_s-K^n_t>0\}\wedge 1$$ where
$K^n_t=\int_{0}^{t}n\left( Y_{s}^{n}-S_{s}\right) ^{-}ds$. Let us
show that
$1_{[\tau^*_t<1]}Y_{\tau^*_t}^{n}=1_{[\tau^*_t<1]}Y_{\tau^*_t}^{n}\wedge
S_{\tau^*_t}^{n}.$
\medskip

Let $\omega$ be fixed such that $\tau^*_t(\omega)<1$. Then there
exists a sequence $(t_k)_{k\geq 0}$ of real numbers which decreases
to $\tau^*_t(\omega)$ such that $Y^n_{t_k}(\omega)\leq
S_{t_k}(\omega)$. As $Y^n$ and $S$ are RCLL processes then taking
the limit as $k\rightarrow \infty$ we obtain $Y^n_{\tau^*_t}\leq
S_{\tau^*_t}$ which implies
$1_{[\tau^*_t<1]}Y_{\tau^*_t}^{n}=1_{[\tau^*_t<1]}Y_{\tau^*_t}^{n}\wedge
S_{\tau^*_t}^{n}.$
\medskip

Now from (\ref{eqstep4}) we deduce that:
$$
\begin{array}{ll}
Y_{t}^{n}&=Y_{\tau^*_t}^{n}+\int_{t}^{\tau^*_t}g(s)ds
-\int_{t}^{\tau^*_t}Z_{s}^{n}dB_{s}
-\int_{t}^{\tau^*_t}\!\!\int_{U}V_{s}^{n}(e) \tilde{\mu
}( ds,de)\\
{}& =1_{[\tau^*_t<1]}Y_{\tau^*_t}^{n}\wedge S_{\tau^*_t}+\xi
1_{[\tau^*_t=1]}+\int_{t}^{\tau^*_t}g(s)ds
-\int_{t}^{\tau^*_t}Z_{s}^{n}dB_{s}
-\int_{t}^{\tau^*_t}\!\!\int_{U}V_{s}^{n}( e) \tilde{\mu }( ds,de) .
\end{array}
$$
Taking the conditional expectation and using inequality
(\ref{eqyn1}) we obtain: $\forall n\geq 0$ and $t\leq 1$,
\begin{equation} \label{eqyn}
Y_t^n=\mbox{ esssup}_{\tau \geq t}\E[\int_t^\tau g(s)ds +(S_\tau
\wedge Y^n_\tau) 1_{[\tau<1]}+\xi 1_{[\tau=1]}|{\cal
F}_t].\Box\end{equation}

\ni\bf{Step 3:} There exists a RCLL process $(Y_t)_{t\leq 1}$ of
${\cal S}^2$ such that: P-$a.s.$,

$(i)$ $Y={\cal H}^{2,1}-lim_{n\rightarrow \infty}Y^n$, $S\leq Y$ and
finally for any $t\leq 1$, $Y^n_t\nearrow Y_t$.

$(ii)$ for any $t\leq 1$, \begin{equation} \label{eqy}Y_t=\mbox{
esssup}_{\tau \geq t}\E[\int_t^\tau g(s)ds +S_\tau 1_{[\tau<1]}+\xi
1_{[\tau=1]}|{\cal F}_t]\;\,\, (Y_1=\xi).\end{equation}

\noindent Actually for $t\leq 1$ let us set
$$\tilde{Y}_t=\mbox{ esssup}_{\tau \geq t}\E[\int_t^\tau g(s)ds +
S_\tau 1_{[\tau<1]}+\xi 1_{[\tau=1]}|{\cal F}_t].$$ The process
$\tilde{Y}$ belongs to ${\cal S}^2$ since $S$ is so, $g\in {\cal
H}^{2,1}$ and $\xi$ is square integrable. On the other hand for any
$n\geq 0$ and $t\leq 1$ we have $Y^n_t\leq \tilde{Y}_t$. Thus there
exits a $\cal P$-measurable process $Y$ such that P-$a.s.$ for any
$t\leq 1$, $Y^n_t\nearrow Y_t\leq \tilde{Y}_t$ and then $Y={\cal
H}^{2,1}-lim_{n\rightarrow \infty}Y^n$. Besides the process
$(Y^n_t+\int_0^tg_sds)_{t\leq 1}$ is a RCLL supermartingale as a
Snell envelope of $(\int_0^t g_sds +(S_t \wedge Y^n_t) 1_{[t<1]}+\xi
1_{[t=1]})_{t\leq 1}$ and it converges increasingly to
$(Y_t+\int_0^tg_sds)_{t\leq 1}$. It follows that this latter process
is also an RCLL supermartingale (see e.g. \cite{dm}, pp.86).
Henceforth the process $Y$ is also RCLL and belongs to ${\cal S}^2$
since it is dominated by $\tilde{Y}$ which is an element of ${\cal
S}^2$.

Next let us prove that $Y\geq S$. Through (\ref{esti1}) we have:
$$
\E[Y_{0}^{n}]=\E[\xi +\int_{0}^{1}g(s)ds]+ \E[\int_{0}^{1}n\left(
Y_{s}^{n}-S_{s}\right) ^{-}ds].$$ Dividing the two hand-sides by $n$
and taking the limit as $n\rightarrow \infty$ to obtain
$E[\int_{0}^{1}( Y_{s}-S_{s})^{-}ds]=0$. As the processes $Y$ and
$S$ are RCLL then $P$-a.s., $Y_t\geq S_t$ for $t<1$. But
$Y_1=\xi\geq S_1$, therefore $Y\geq S$.

Finally let us show that $Y$ satisfies (\ref{eqy}). But this a
direct consequence of the continuity of the Snell envelope through
sequences of increasing RCLL processes (see Appendix [A1]). Actually
on the one hand, the sequence of increasing RCLL processes
$(S_t\wedge Y^n_t)1_{[t<1]}+\xi 1_{[t=1]})_{t\leq 1}$ converges
increasingly to the RCLL $(S_t1_{[t<1]}+\xi 1_{[t=1]})_{t\leq 1}$
since $Y\geq S$. Therefore because of (\ref{eqyn}) the sequence
$(\int_0^t g_sds +Y^n_t)_{t\leq 1}$ converges to $ \mbox{
esssup}_{\tau \geq t}\E[\int_0^\tau g(s)ds + S_\tau 1_{[\tau<1]}+\xi
1_{[\tau=1]}|{\cal F}_t]$ which then is equal to $(\int_0^t g_sds
+Y_t)_{t\leq 1}$ and which implies that $Y$ satisfies (\ref{eqy}).
$\Box$ \bigskip

\ni \bf{Step 4}: There exist three processes $Z\in {\cal H}^{2,d}$,
$V\in {\cal L}^2$ and $K\in {\cal S}^2_i$ such that $(Y,Z,V,K)$ is a
the solution of the BSDE associated with $(g,\xi,S)$.  \bigskip

We know from (\ref{eqy}), that the process
$(\int_0^tg_sds+Y_t)_{t\leq 1}$ is a Snell envelope. Then through
Appendix [A2], there exit a process $K\in {\cal S}_i^2$ ($K_0=0$)
and an ${\cal F}_t$-martingale $(M_t)_{t\leq 1}$ which belongs to
${\cal S}^2$ such that:
$$\forall t\leq 1,\,\,\int_0^tg_sds+Y_t=M_t-K_t.$$
Additionally $K=K^c+K^d$ where $K^c$ is continuous non-decreasing
and $K^d$ non-decreasing purely discontinuous predictable and such
that for any $t\leq 1$, $\Delta_t
K^d=(S_{t-}-Y_t)^+1_{[Y_{t-}=S_{t-}]}$.

Now the martingale $M$ belongs to ${\cal S}^2$ then the
representation property (see e.g. \cite{iw}) implies the existence
of two processes $Z$ and $V$ which belong respectively to ${\cal
H}^{2,d}$ and ${\cal L}^2$ such that:
$$P-a.s., \,\,\forall t\leq 1,\,\,
M_t=M_0+\int_0^t\{Z_sdB_s+\int_UV_s(e)\tilde{\mu}(ds,de)\}.$$

Let us now show that $\int_0^1(Y_s-S_s)dK^c_s=0$. First let us
remark that the Snell envelope of
$(\int_0^tg_sds+S_t1_{[t<1]}+\xi1_{[t=1]}+K^d_t)_{t\leq 1}$ is
nothing else but $(\int_0^tg_sds+Y_t+K^d_t)_{t\leq 1}$.

Actually for any $t\leq 1$ we have
$\int_0^tg_sds+Y_t=M_t-K^c_t-K^d_t$, therefore the process
$(\int_0^tg_sds+Y_t+K^d_t)_{t\leq 1}$ is also a RCLL supermartingale
which dominates the process
$(\int_0^tg_sds+S_t1_{[t<1]}+\xi1_{[t=1]}+K^d_t)_{t\leq 1}$. Besides
if $(N_t)_{t\leq T}$ is a supermartingale of class [D] which
dominates this latter process then $(N_t-K^d_t)_{t\leq T}$ still a
supermartingale of class [D] which is greater than
$(\int_0^tg_sds+S_t1_{[t<1]}+\xi1_{[t=1]})_{t\leq 1}$. Therefore
$P-a.s.$ for any $t\leq 1$ we have $N_t-K^d_t\geq \int_0^tg_sds+Y_t$
which implies that $\forall t\leq 1$, $N_t\geq
\int_0^tg_sds+Y_t+K^d_t$. It means that the process
$(\int_0^tg_sds+Y_t+K^d_t)_{t\leq T}$ is the smallest
supermartingale of class [D] which dominates
$(\int_0^tg_sds+S_t1_{[t<1]}+\xi1_{[t=1]}+K^d_t)_{t\leq 1}$ and then
it is its Snell envelope.
\medskip

Next the Snell envelope $(\int_0^tg_sds+Y_t+K^d_t=M_t-K^c_t)_{t\leq
1}$ of the process
$(\int_0^tg_sds+S_t1_{[t<1]}+\xi1_{[t=1]}+K^d_t)_{t\leq 1}$ is
regular (see Appendix A3) then for any $t\leq 1$, the stopping time
$\tau_t=\inf\{s\geq t, K_s>K_t\}\wedge 1$ is optimal (see [A3])
therefore we have
$\int_t^{\tau_t}(Y_s+K^d_s-S_s-K^d_s)dK^c_s=\int_t^{\tau_t}(Y_s-S_s)dK_s^c=0$.
As $t$ is arbitrary then we have $\int_0^{1}(Y_s-S_s)dK_s^c=0$.
\medskip

Collecting now all those properties yields that the quadruple
$(Y,Z,V,K)$ is a solution for the BSDE associated with $(g,\xi,S)$,
i.e., \begin{equation}\label{eqg}\left\{
\begin{array}{ll}
{}& Y\in {\cal S}^{2},\,Z\in {\cal H}^{2,d},\,V\in {\cal
L}^{2}\mbox{ and
}K\in {\cal S}^{2}_i \\
{}& Y_{t}=\xi +\displaystyle \int_{t}^{1}g(s)ds+K_{1}-K_{t}-
\int_{t}^{1}Z_{s}dB_{s}-\displaystyle \int_{t}^{1}\!\!\!\!
\int_{U}^{{}}V_{s}(e)\tilde{\mu}(ds,de)\,,\forall t\leq 1 \\
{}& Y\geq S \mbox{ and if $K^c$ (resp. $K^d$) is the continuous
(resp. purely discontinuous) part of $K$,}\\& \mbox{then   $K^d$ is
${\cal P}^d$-measurable, }
\int_{0}^{1}(Y_{t}-S_{t})dK_{t}^{c}=0\,\,\text{and }\forall t\leq 1,
\Delta K_t^d=(S_{t-}-Y_t)^{+}1_{[ Y_{t-}=S_{t-}]}. \,\Box
\end{array}
\right.
\end{equation}

\begin{remark}: (i) In using the generalization of the monotonic
limit theorem of S.Peng (\cite{monlim}), one can show that the
sequence $(Z^n)_{n\geq 0}$ (resp. $(V^n)_{n\geq 0}$) converge weakly
in ${\cal H}^{2,d}$ (resp. $\text{L}^2(dt\otimes dP\otimes
d\lambda)$) and strongly in $\text{L}^p(dt\otimes dP)$ (resp.
$\text{L}^2(dt\otimes dP\otimes d\lambda)$) for any $p<2$. On the
other hand for any stopping time $\tau$ the sequence
$(K^n_\tau)_{n\geq 0}$ converges in $L^p(dP)$ to $K_\tau$.
\medskip

(ii) It can be easily seen that we would have been able to construct
the solution of the BSDE (\ref{eqg}) directly from (\ref{eqy}) in
defining the process $Y$, and then $Z,\,V,\, K$ as in Step 4 above.
This procedure is the well known Snell envelope method. However we
think that proving the convergence of the penalization scheme
(\ref{esti1}) could be helpful especially when dealing with, either
numerical schemes for BSDEs of type (\ref{eqg}) or viscosity
solutions of PDIEs with discontinuous obstacles. $\Box$
\end{remark}

We are now ready to give the main result of this section. \medskip

\begin{theorem}: The reflected BSDE with generalized jumps (\ref{prince})
associated with $(f,\xi,S)$ has a unique solution $(Y,Z,K,V)$.
\end{theorem}
$Proof$: It remains to show existence which will be obtained via a
fixed point argument. Actually let ${\cal D}:={\cal H}^{2,1}\times
{\cal H}^{2,d}\times {\cal L}^{2}$ endowed with the norm
\[
\Vert (Y,Z,V)\Vert _{\alpha }=\{I\!\!E[\displaystyle
\int_{0}^{1}e^{\alpha s}(|Y_{s}|^{2}+|Z_{s}|^{2}+\displaystyle
\int_{U}^{{}}|V_{s}(e)|^{2}\lambda (de))ds]\}^{1/2};\,\,\alpha >0.
\] On the other hand let $\Phi $ be the map from ${\cal D}$ into itself which
with $(Y,Z,V)$ associates \\ $\Phi
(Y,Z,V)=(\tilde{Y},\tilde{Z},\tilde{V})$ where
$(\tilde{Y},%
\tilde{Z},\tilde{K},\tilde{V})$ is the solution of the reflected
BSDE associated with $(f(t,Y_{t},Z_{t},V_{t}),\xi ,S)$. Let
$(Y^{\prime },Z^{\prime },V^{\prime })$ be another triple of ${\cal D}$ and $%
\Phi (Y^{\prime },Z^{\prime },V^{\prime })=(\tilde{Y}^{\prime },\tilde{Z}%
^{\prime },\tilde{V}^{\prime })$. Using It\^{o}'s formula we obtain:
$\forall t\leq 1$,
\[
\begin{array}{ll}
& e^{\alpha t}(\tilde{Y}_{t}-\tilde{Y}_{t}^{\prime })^{2}+\alpha %
\displaystyle \int_{t}^{1}e^{\alpha
s}(\tilde{Y}_{s}-\tilde{Y}_{s}^{\prime
})^{2}ds+\displaystyle \int_{t}^{1}e^{\alpha s}|\tilde{Z}_{s}-\tilde{Z}%
_{s}^{\prime }|^{2}ds+ \\
& \qquad \sum_{t<s\leq 1}e^{\alpha s}(\Delta _{s}\tilde{Y}-\Delta
_{s}\tilde{Y}%
^{\prime })^{2}=(M_{1}-M_{t})+2\displaystyle \int_{t}^{1}e^{\alpha s}(\tilde{%
Y}_{s-}-\tilde{Y}_{s-}^{\prime })(d\tilde{K}_{s}-d\tilde{K}_{s}^{\prime }) \\
& \qquad \qquad +2\displaystyle \int_{t}^{1}e^{\alpha s}(\tilde{Y}_{s}-%
\tilde{Y}_{s}^{\prime })(f(s,Y_{s},Z_{s},V_{s})-f(s,Y_{s}^{\prime
},Z_{s}^{\prime },V_{s}^{\prime }))ds%
\end{array}
\]
where $(M_{t})_{t\leq 1}$ is a martingale. But for any $t\leq 1$,
$\displaystyle \int_{t}^{1}e^{\alpha
s}(\tilde{Y}_{s-}-\tilde{Y}_{s-}^{\prime
})(d\tilde{K}%
_{s}-d\tilde{K}_{s}^{\prime })\leq 0. $ This can be shown as in the
proof of uniqueness in Proposition \ref{unicite}. Therefore taking
expectation in both hand-sides yields
\[
\begin{array}{ll}
& \alpha I\!\!E[\displaystyle \int_{t}^{1}e^{\alpha s}(\tilde{Y}_{s}-\tilde{Y%
}_{s}^{\prime })^{2}ds]+I\!\!E[\displaystyle \int_{t}^{1}e^{\alpha s}|\tilde{%
Z}_{s}-\tilde{Z}_{s}^{\prime }|^{2}ds]+I\!\!E[\displaystyle %
\int_{t}^{1}e^{\alpha s}ds\displaystyle \int_{U}^{{}}(\tilde{V}_{s}(e)-%
\tilde{V}_{s}^{\prime }(e))^{2}\lambda (de)] \\
& \qquad\leq 2I\!\!E[\displaystyle \int_{t}^{1}e^{\alpha
s}(\tilde{Y}_{s}-\tilde{Y}%
_{s}^{\prime })(f(s,Y_{s},Z_{s},V_{s})-f(s,Y_{s}^{\prime
},Z_{s}^{\prime
},V_{s}^{\prime }))ds] \\
& \qquad\leq C_f\epsilon I\!\!E[\displaystyle \int_{t}^{1}e^{\alpha
s}(\tilde{Y}_{s}-%
\tilde{Y}_{s}^{\prime })^{2}ds]+\frac{C_f}{\epsilon }I\!\!E[\displaystyle %
\int_{t}^{1}e^{\alpha s}\{|Y_{s}-Y_{s}^{\prime
}|^{2}+|Z_{s}-Z_{s}^{\prime
}|^{2}+ \\
& \qquad \qquad \qquad \qquad \qquad \displaystyle %
\int_{U}^{{}}|V_{s}(e)-V_{s}^{\prime }(e)|^{2}\lambda (de)\}ds].%
\end{array}
\]
It implies that
\[
\begin{array}{l}
(\alpha -\epsilon C_f)I\!\!E[\displaystyle \int_{t}^{1}e^{\alpha
s}(\tilde{Y} _{s}-\tilde{Y}_{s}^{\prime })^{2}ds]+
I\!\!E[\displaystyle \int_{t}^{1}e^{\alpha
s}|\tilde{Z}_{s}-\tilde{Z}_{s}^{\prime }|^{2}ds]+\E[\displaystyle
\int_{t}^{1}e^{\alpha s}ds \displaystyle
\int_{U}^{{}}(\tilde{V}_{s}(e)-\tilde{V}_{s}^{\prime
}(e))^{2}\lambda (de)] \\\qquad \qquad \leq \frac{C_f}{\epsilon
}I\!\!E[\displaystyle \int_{t}^{1}e^{ \alpha
s}\{|Y_{s}-Y_{s}^{\prime }|^{2}+|Z_{s}-Z_{s}^{\prime }|^{2}+
\displaystyle \int_{U}^{{}}|V_{s}(e)-V_{s}^{\prime }(e)|^{2}\lambda
(de)\}ds].
\end{array}
\]
Now let $\alpha $ be great enough and $\epsilon $ such that
$C_f<\epsilon
<\frac{%
\alpha -1}{C_f}$, then $\Phi $ is a contraction on ${\cal D}$.
Therefore there exists a triple $(Y,Z,V)$ such that $\Phi
(Y,Z,V)=(Y,Z,V)$ and then with  $K$ the quadruple $(Y,Z,V,K)$ is the
unique solution of the reflected BSDE associated with $(f,\xi ,S)$
since $Y\in {\cal S}^2$. $\Box $

\subsection{On the regularity of the process $K$}

We now focus on the regularity of the process $K$ of the solution of
the BSDE (\ref{prince}).
\begin{proposition}
Let $\left( Y,Z,K,V\right) $ be the unique solution of the reflected
BSDE (\ref{eqg}) associated with $(g(t),\xi,S)$ and let us consider
the following assertions:

$(i)$ $K$ is continuous $(K^{d}=0)$

$(ii)$ the sequence of processes $(Y^{n})_{n\geq 0}$ of
(\ref{esti1}) converges uniformly to $Y$

$(iii)$ $I\!\!E[ \int_{0}^{1}\left| Z_{s}^{n}-Z_{s}\right|
^{2}ds+\int_{0}^{1}\int_{U}\left| V_{s}^{n}\left( e\right)
-V_{s}\left( e\right) \right| ^{2}\lambda \left( de\right) ds]
\rightarrow 0.$

\noindent Then it holds true that $(i)$ and $(ii)$ are equivalent
and the statement $(ii)$ implies $(iii)$.
\end{proposition}
$Proof$:  $(i)\Rightarrow (ii)$: Let us assume that $(i)$ is
fulfilled. Therefore the jumps of $Y$ are the same as the ones of
its Poisson martingale part. It follows that for any $t\leq 1$,
$^{p}Y_t=Y_{t-}$ where $^pY$ is the predictable projection of $Y$.
Because in that case $Y$ has only inaccessible jumps. But
$Y^{n}\nearrow Y$ , thus $^{p}Y^{n}\nearrow {^pY}$, then for any
$t\leq 1$ we have $Y^n_{t-}\nearrow Y_{t-}$. It follows from a
generalized Dini theorem (see \cite{dm}, pp.203) that P-a.s. the
sequence $(Y^n)_{n\geq 0}$ converges uniformly to $Y$.
\medskip

\noindent $(ii)\Rightarrow (i)$: If $(Y^{n})_{n\geq 0}$ converges
uniformly to $Y$ then $ Y_{-}^{n}\nearrow Y_{-}$ and
$^{p}Y^{n}\nearrow $ $^{p}Y$. But $Y^{n}$ is RCLL and has only
inaccessible jumps, then $^{p}Y^{n}=Y_{-}^{n}$. It follows that
$^{p}Y=Y_{-}$. Therefore the Snell envelope
$(\int_0^tf(s,Y_s,Z_s,V_s)ds+Y_t)_{t\leq 1}$ is regular, thus
$K^d=0$ (see Appendix [A3]) which yields the desired result.
\medskip

\noindent $(ii)\Rightarrow (iii)$: If the sequence $(Y^n)_{n\geq 0}$
converges uniformly to $Y$ then through dominated convergence
theorem we have also $\E[\sup_{t\leq 1}(Y-Y^{n})^{2}]\searrow 0$. So
to obtain the result, it is enough to apply It\^o's formula with
$(Y-Y^n)^2$ after having remarked that $\E[\sup_{t\leq
1}((S_t-Y^{n})^+)^{2}]\rightarrow 0$ since $Y\geq S$. Therefore we
have:
\[
\E[\int_{0}^{1}|Z_{s}^{n}-Z_{s}|^{2}ds+\int_{0}^{1}\!\!\!\int_{U}|
V_{s}^{n}(e) -V_{s}(e)|^{2}\lambda (de)ds] \rightarrow 0. \Box
\]
\begin{remark}: let us point out that the implication $(iii)\Rightarrow (i)$
is not true in general. Actually let us consider the following
deterministic counter-example. Assume $\xi=\frac{1}{2}$, $f\equiv 0$
and $S_t=1_{[t<\frac{1}{2}]}$. Therefore the solution of the BSDE
associated with $(f,\xi,S)$ is $Z\equiv 0$, $V\equiv 0$,
$Y_t=1_{[t<\frac{1}{2}]}+\frac{1}{2}1_{[t\geq \frac{1}{2}]}$ and
$K_t=K^d_t=\frac{1}{2}1_{[\frac{1}{2}\leq t\leq 1]}$. Indeed it is
easily seen that for any $t\leq 1$, we have
$Y_{t}=\frac{1}{2}+K_{1}^{d}-K_{t}^{d}$. So obviously the statement
$(iii)$ holds true since $Z^n=Z=0$, $V^n=V=0$ but $K$ is not
continuous. $\Box$
\end{remark}

\section{BSDEs with two discontinuous reflecting barriers}

We now consider the problem of reflection with respect to two
barriers, an upper and a lower ones. So let us give two processes
$L:=(L_{t})_{t\leq T}$ and $U:=(U_{t})_{t\leq }$ which stand for the
barriers where "the solution" is reflected and which satisfy:
\medskip

$(i)$ $L$ and $U$ belong to ${\cal S}^2$ and $P-a.s.,\,\,\forall
t\leq 1, \,L_t\leq U_t$ and $L_1\leq \xi\leq U_1$

$(ii)$ there exist two non-negative supermartingales
$(h_t)_{t\leq1}$ and $(h'_t)_{t\leq 1}$ of ${\cal S}^2$ such that:
\[ \forall t\leq 1,\,\,L_{t}\leq
h_t-h'_t\leq U_{t}.
\]This assumption is the so-called \it{Mokobodski's hypothesis}

$(iii)$ $\forall t< 1, L_{t-}<U_{t-}$ and $L_t<U_t$. $\Box$
\bigskip

Let us now introduce the BSDE associated with $(f,\xi,L,U)$. A
solution is a quintuple of processes
$(Y_{t},Z_{t},K_{t}^{+},K_{t}^{-},V_{t})_{t\leq 1}$ which satisfies:
\begin{equation}
\label{2-b} \left\{
\begin{array}{l}
(i)\,Y\in {\cal S}^{2},K^{\pm }\in {\cal S}^2_i,Z\in {\cal
H}^{2,d}\mbox{ and
}V\in {\cal L}^2 \\
(ii)\,-dY_{t}=f(t,Y_{t},Z_{t},V_t)dt+dK_{t}^{+}-dK_{t}^{-}-Z_{t}dB_{t}-
\displaystyle \int_{U}\!\!V_{t}(e)\tilde{\mu}(dt,de),\,t\leq
1;\,\,Y_{1}=\xi
\\
(iii)\,\forall t\leq 1,L_{t}\leq Y_{t}\leq U_{t}\mbox{ and if
}K^{\pm c} \mbox{ is the continuous part of }K^{\pm }\mbox{ then
}(Y_{t}-L_{t})dK_{t}^{+c}=0 \\ \qquad\mbox{ and
}(U_{t}-Y_{t})dK_{t}^{-c}=0 \\
(iv)\,\mbox{if $K^{\pm d}$ denotes the purely discontinuous part of
$K^\pm$ then $K^{\pm d}$ is ${\cal P}^d$-measurable }\\ \mbox{ and
$\forall t\leq 1, \Delta K_t^{+d}=(L_{t-}-Y_t)^{+}1_{[
Y_{t-}=L_{t-}]}$ and $\Delta K_t^{-d}=(Y_t-U_{t-})^{+}1_{[
Y_{t-}=U_{t-}]}$}.
\end{array}
\right.
\end{equation}
Note that the BSDE (\ref{2-b}) may have not a solution. Actually if
for example $L$ is not a semimartingale, and $L$, $U$ coincide then
obviously the equation cannot have a solution since $Y$ is a
semimartingale. $\Box$
\medskip

First we are going to focus on the issue of uniqueness of the
solution of (\ref{2-b}).
\begin{proposition}: \bf{Uniqueness}
\medskip

\noindent The BSDE with two reflecting barriers associated with
$(f,\xi,L,U)$  (\ref{2-b}) has at most one
solution.\end{proposition}

\noindent $Proof$: Assume that $(Y,Z,K^\pm,V)$ and $(Y^{\prime
},Z^{\prime },K^{\prime \pm },V^{\prime })$ are two solutions of
(\ref{2-b}). First let us show that for any $t\leq 1$,
$\int_t^1(Y_{s-}-Y_{s-}^{\prime }{})(dK_{s}-dK_{s}^{\prime })\leq 0$
where $K=K^+-K^-$ and $K'={K'}^+-{K'}^-$.

For any $t\leq 1$ we have,
$$
\int_{]t,1]}^{{}}(Y_{s-}-Y_{s-}^{\prime }{})(dK_{s}-dK_{s}^{\prime
})=
\int_{]t,1]}^{{}}(Y_{s-}-Y_{s-}^{\prime})(dK_{s}^{c}-dK_s^{\prime
c})+ \int_{]t,1]}^{{}}(Y_{s-}-Y_{s-}^{\prime
}{})(dK_{s}^{d}-dK_{s}^{\prime d }).
$$
The processes $Y$ and $Y^{\prime }$ belong to ${\cal S}^{2}$ and
their jumps
$%
\delta (\omega ):=\{t\in [0,1],\Delta _{t}Y(\omega)\neq 0\}$ and
$\delta ^{\prime }(\omega ):=\{t\in [0,1],\Delta _{t}Y^{\prime
}(\omega)\neq 0\}$ are at most countable. Therefore
\[
\begin{array}{ll}
\!\!\!\displaystyle \int_{]t,1]}^{{}}(Y_{s-}-Y_{s-}^{\prime
})(dK_{s}^{c}-dK_s^{\prime c})=\!\!\!\displaystyle
\int_{]t,1]}^{{}}(Y_{s}-Y_{s}^{\prime })(dK_{s}^{+c}-dK_s^{\prime +
c })- \!\!\!\displaystyle \int_{]t,1]}^{{}}(Y_{s}-Y_{s}^{\prime
})(dK_{s}^{-c}-dK_s^{\prime -c}).
\end{array}
\]
But since for any $t\leq 1$, $(Y_t-L_t)dK^{+c}_t= 0$ then
$$\begin{array}{ll}\displaystyle
\int_{]t,1]}^{{}}(Y_{s}-Y_{s}^{\prime})(dK_{s}^{+c}-dK_s^{\prime
+c})=-\displaystyle \int_{]t,1]}^{}(Y_{s}-L_{s})dK_{s}^{\prime
+c}+\displaystyle \int_{]t,1]}^{}(L_{s}-Y_{s}^{\prime
})d{K}_{s}^{+c}\leq 0.
\end{array}
$$
In the same way, we can also show that $
\int_{]t,1]}^{{}}(Y_{s}-Y_{s}^{\prime })(dK_{s}^{-c}-dK_s^{\prime
-c})\geq 0$. Therefore we obtain:
\begin{equation}
\label{inedb} \forall t\leq
1,\,\,\begin{array}{c}\int_{]t,1]}^{{}}(Y_{s-}-Y_{s-}^{\prime
})(dK_{s}^{c}-dK_s^{\prime c})\leq 0.\end{array}
\end{equation}
Let us now focus on the discontinuous parts of $K-K'$. For any
$t\leq 1$,
\begin{equation}
\label{inebc}
\begin{array}{ll}
\displaystyle \int_{]t,1]}^{{}}(Y_{s-}-Y_{s-}^{\prime
})(dK_{s}^{d}-dK_s^{\prime d})&=\displaystyle
\int_{]t,1]}^{{}}(Y_{s-}-Y_{s-}^{\prime})(dK_{s}^{+d}-dK_s^{\prime +
d})-\\{}&\qquad\qquad \displaystyle
\int_{]t,1]}^{{}}(Y_{s-}-Y_{s-}^{\prime })(dK_{s}^{-d}-dK_s^{\prime
-d}).\end{array}
\end{equation}
But
\[
\int_{]t,1]}^{{}}(Y_{s-}-Y'_{s-})(dK_{s}^{+d}-dK_s^{\prime +d})=
\int_{]t,1]}(L_{s-}-Y'_{s-})dK_{s}^{+d}
-\int_{]t,1]}^{{}}(Y_{s-}-L_{s-})dK_s^{\prime +d} \leq 0
\]
since $Y\geq L$ (resp. $Y'\geq L$) and the jumps of $K^{+d}$ (resp.
${K'}^{+d}$) occur only when $Y_-=L_-$ (resp. $Y'_-=L_-$). In the
same way we have $\int_{]t,1]}^{{}}(Y_{s-}-Y_{s-}^{\prime
})(dK_{s}^{-d}-dK_s^{\prime -d})\geq 0$. It follows from
(\ref{inebc}) that $\int_{]t,1]}^{{}}(Y_{s-}-Y_{s-}^{\prime
})(dK_{s}^{d}-dK_s^{\prime d})\leq 0.$ Combining now this inequality
with (\ref{inedb}) we deduce that for any $t\leq 1$, we have
$\int_{]t,1]}(Y_{s-}-Y_{s-}^{\prime })(dK_{s}-dK_s^{\prime })\leq
0$.
\medskip

Now using It\^{o}'s formula with $(Y-Y^{\prime })^2$ and following
the same steps as in the proof of uniqueness of BSDEs with one
reflecting barrier (see Proposition \ref{unicite}) to obtain that
$Y=Y'$, $Z=Z'$, $V=V'$ and finally $K=K'$. Thus, due to their
expressions, we have also $K^{+d}={K'}^{+d}$ and $K^{-d}={K'}^{-d}$
and then $K^{+c}-K^{-c}={K'}^{+c}-{K'}^{-c}$. It remains to show
that $K^{+c}={K'}^{+c}$ and $K^{-c}={K'}^{-c}$.

Indeed for any $t\leq 1$ we have:
$$
\begin{array}{ll}
\int_{0}^{t}(Y_s-L_s)d(K^{+c}_s-K^{-c}_s)&=-\int_{0}^{t}(U_s-L_s)dK_s^{-c}=
\int_{0}^{t}(Y'_s-L_s)d(K_s^{\prime +c}-K_s^{\prime
-c})\\{}&=-\int_{0}^{t}(U_s-L_s)dK_s^{\prime -c}.
\end{array}
$$
Therefore $K^{-c}={K'}^{-c}$ since for any $t<T$, $L_t<U_t$ and then
$K^{+c}={K'}^{+c}$. Thus we have uniqueness of the solution. $\Box$
\medskip

Once again to show that equation (\ref{2-b}) has a solution we are
going first to suppose that $f$ does not depend on $(y,z,v)$, i.e,
$f(t,\omega,y,z,v)=f(t)$. Then we have the following:

\begin{theorem}:
There exists a unique 5-uple of processes
$(Y_{t},Z_{t},K_{t}^{+},K_{t}^{-},V_{t})_{t\leq T}$ solution of the
backward stochastic differential equation with two reflecting
barriers associated with $(f(t),\xi,L,U)$.
\end{theorem}
$Proof$: Even if the barriers have predictable jumps, the proof of
this theorem, in its main steps, is classical (see e.g. \cite{CK},
\cite{hl}).

Let us consider the following processes defined by: $\forall t\leq
1$,
$$
\begin{array}{c}
H_t=(h_t+\E[\xi^-|F_t])1_{[t<1]}+\E[\int_{t}^{1}f(s)^-ds|{\cal F}_t],\\
\Theta_t=(h'_t+\E[\xi^+|F_t])1_{[t<1]}+\E[\int_{t}^{1}f(s)^+ds|{\cal F}_t],\\
\tilde{L}_t=L_t1_{[t<1]}+\xi1_{[t=1]}-\E[\xi+\int_{t}^{1}f(s)ds|{\cal
F}_t]
\\
\mbox{ and
}\tilde{U}_t=U_t1_{[t<1]}+\xi1_{[t=1]}-\E[\xi+\int_{t}^{1}f(s)ds|{\cal
F}_t],
\end{array}
$$ where $f(t)^-=\max\{-f(t),0\}$ and $f(t)^+=\max\{f(t),0\}$. Since $h$
and $h'$ are non-negative supermartingales then $H$ and $\Theta$ are
also non-negative supermartingales which moreover belong to ${\cal
S}^2$ and verify $H_1=\Theta_1=0$. On the other hand, through
Mokobodski's hypothesis, we can easily verify that for any $t\leq 1$
we have:
\begin{equation}\label {eqmoko}
\tilde{L}\leq H-\Theta\leq \tilde{U}.
\end{equation}
Next let us consider the sequences $(N^\pm_n)_{n\geq 0}$ of
processes defined recursively as follows: $$ N_0^\pm=0 \mbox{ and
for }n\geq 0, N^{+,n+1}=R(N^{-,n}+\tilde{L}) \mbox{ and
}N^{-,n+1}=R(N^{+,n}-\tilde{U})$$ where $R$ is the Snell envelope
operator (see Appendix). Now by induction and using (\ref{eqmoko})
we can easily verify that:
$$ \forall n\geq 0,\,\, 0\leq N^{+,n}\leq N^{+,n+1}\leq H\mbox{ and
} 0\leq N^{-,n}\leq N^{-,n+1}\leq \Theta.$$ It follows that the
sequence $(N^+_n)_{n\geq 0}$ (resp. $(N^-_n)_{n\geq 0}$) converges
pointwisely to a supermartingale $N^+$ (resp. $N^-$) (see e.g.
\cite{dm}, pp.86). In addition $N^+$ and $N^-$ belong to ${\cal
S}^2$ and verify (see [A1]) : $$ N^+=R( N^-+\tilde{L}) \mbox{ and
}N^{-}=R(N^{+}-\tilde{U}).$$ Next the Doob-Meyer decompositions of
$N^\pm$ yield : $$ \forall t\leq 1, N^\pm_t=M^\pm_t-K^\pm_t$$ where
$M^\pm$ are RCLL martingales and $K^\pm$ non-decreasing processes
such that $K^\pm_0=0$. Moreover since $N^\pm\in {\cal S}^2$ then
$E[(K^\pm_T)^2]<\infty$ (see [A2]). Therefore $M^\pm$ belong also to
${\cal S}^2$ and then there exist processes $Z^\pm\in {\cal
H}^{2,d}$ and $V^\pm \in {\cal L}^2$ such that (see \cite{iw}):
$$ \forall t\leq 1,
M^\pm_t=M^\pm_0+\integ{0}{t}\{Z^\pm_sdB_s+\integ{U}{}V^\pm_s(e)\tilde{\mu}(ds,de)\}.$$
Next let us denote by $K^{\pm d}$ (resp. $K^{\pm c}$) the purely
discontinuous (resp. continuous) part of $K^\pm$. In the same way as
shown for BSDEs with one reflecting barrier (see Section 4, Step 4)
we have:
\begin{equation}
\int_0^1(N^+- N^--\tilde{L})dK^{+c}_s=\int_0^1(N^--
N^++\tilde{L})dK^{-c}_s=0.\end{equation} On the other hand the
processes $K^{\pm d}$ are predictable and if $\tau$ is a predictable
stopping time (see [A2]) then $$\{\Delta K^{+d}_\tau>0\}\subset
\{N^+_{\tau -}=N^-_{\tau -}+\tilde{L}_{\tau -}\}\mbox{ and }
\{\Delta K^{-d}_\tau>0\}\subset \{N^-_{\tau -}=N^+_{\tau
-}-\tilde{U}_{\tau -}\}.$$ But $L_-<U_-$ and $\tau$ is predictable
then we have $\tilde{L}_{\tau -}<\tilde{U}_{\tau -}$ since, the
jumps of martingales with respect to $({\cal F})_{t\leq 1}$ are
inaccessible because they come only  from the Poisson part.
Therefore the predictable processes $K^{+d}$ and $K^{-d}$ cannot
jump in the same time otherwise we would have $\tilde{L}_{\tau
-}=\tilde{U}_{\tau -}$ which is impossible. Henceforth
$$\{\Delta
K^{+d}_\tau>0\}\subset \{N^+_{\tau -}=N^-_{\tau }+\tilde{L}_{\tau
-}\}\mbox{ and } \{\Delta K^{-d}_\tau>0\}\subset \{N^-_{\tau
-}=N^+_{\tau }-\tilde{U}_{\tau -}\}.$$ It follows that for any
$t\leq 1$ we have
$$\Delta K^{+d}_t=(N^+_{t-}-N^+_{t})^+1_{\{N^+_{t-}=N^-_{t
-}+\tilde{L}_{t -}\}}=(N^-_{t}+\tilde{L}_{t
-}-N^+_{t})^+1_{\{N^+_{t-}=N^-_{t -}+\tilde{L}_{t -}\}}.$$ In the
same way we obtain:
$$\forall t\leq 1,\,\,\Delta K^{-d}_t=(N^-_{t-}-N^-_{t})^+1_{\{N^-_{t-}=N^+_{t
-}-\tilde{U}_{t -}\}}=(N^+_{t}-\tilde{U}_{t
-}-N^-_{t})^+1_{\{N^-_{t-}=N^+_{t -}-\tilde{U}_{t -}\}}.$$ Finally
for $t\leq 1$, let us set: $$
Y_t=N^+_t-N^-_t+\E[\xi+\integ{t}{1}f(s)ds|{\cal F}_t], \,\,
Z_t=Z^+_t-Z^-_t+\eta_t,\,\, V_t=V^+_t-V^-_t+\rho_t$$ where the
processes $\eta$ and $\rho$ are such that $$ \forall t\leq 1,\,\,
\E[\xi+\integ{0}{1}f(s)ds|{\cal F}_t]=\E[\xi+\integ{0}{1}f(s)ds]+
\integ{0}{t}\{\eta_sdB_s+\integ{U}{}\rho_s(e)\tilde{\mu}(ds,de)\}.$$
Therefore the quintuple $(Y,Z,V,K^{+},K^{-})$ is the solution of the
BSDE with two reflecting barriers associated with $(f(t),\xi, L,U)$,
i.e.,
$$
\left\{
\begin{array}{l}
\,Y\in {\cal S}^{2},K^{\pm }\in {\cal S}^2_i,Z\in {\cal
H}^{2,d}\mbox{ and }V\in
{\cal L}^2; \\
\,-dY_{t}=f(t,Y_{t},Z_{t},V_t)dt+dK_{t}^{+}-dK_{t}^{-}-Z_{t}dB_{t}-
\displaystyle \int_{U}\!\!V_{t}(e)\tilde{\mu}(dt,de),\,t\leq
1;\,\,Y_{1}=\xi;
\\
\,\forall t\leq 1,L_{t}\leq Y_{t}\leq U_{t}\mbox{ and if }K^{\pm c}
\mbox{ is the continuous part of }K^{\pm }\mbox{ then
}(Y_{t}-L_{t})dK_{t}^{+c}=0 \\ \qquad\mbox{ and
}(U_{t}-Y_{t})dK_{t}^{-c}=0 ;\\
\,\mbox{if $K^{\pm d}$ denotes the purely discontinuous part of
$K^\pm$ then $K^{\pm d}$ is ${\cal P}^d$-measurable }\\ \mbox{
$\,\,\,$ and $\forall t\leq 1, \Delta K_t^{+d}=(L_{t-}-Y_t)^{+}1_{[
Y_{t-}=L_{t-}]}$  and $\Delta K_t^{-d}=(Y_t-U_{t-})^{+}1_{[
Y_{t-}=U_{t-}]}$}.
\end{array}
\right.
$$\vskip 0.5cm

We are now ready to give the main result of this section.

\begin{theorem}\label{thm3}: The reflected BSDE
(\ref{2-b}) associated with $(f(t,y,z,v), \xi,L,U)$ has a unique
solution $(Y,Z,V,K^+,K^-)$.
\end{theorem}

\ni $Proof$: We give a brief proof since once more it is somehow
classical. Let ${\cal H}:={\cal H}^{2,1}\times {\cal H}^{2,d}\times
{\cal L}^2$ and $\Phi$ be the following application:
$$
\begin{array}{ll}
\Phi:&{\cal H}\longrightarrow {\cal H} \\ (y,z,v):&\mapsto
\Phi(y,z,v)=(\bar{Y},\bar{Z},\bar{V})
\end{array}$$
where $(\bar{Y},\bar{Z},\bar{V})$ is the triple for which there
exists two other processes $\bar{K}^\pm$ which belong to ${\cal
S}^2$ such that $(\bar{Y},\bar{Z},\bar{V},\bar{K}^+,\bar{K}^-)$ is a
solution for the BSDE with two reflecting barriers associated with
$(f(t,y_t,z_t,v_t),\xi,L,U)$. Now let $\alpha >0$, $(y',v',z')\in
{\cal H}$ and $(\bar{Y}',\bar{Z}',\bar{V}')=\Phi (y',z',v')$. Using
It\^o's formula and taking into account that $e^{\alpha
s}(\bar{Y}_s-\bar{Y}'_s)d(\bar{K}^+_s-\bar{K}^-_s-\bar{K'}^+_s+\bar{K'}^-_s)\leq
0$ we show the existence of a constant $\bar{C}<1$, by an
appropriate choice of $\alpha$ (see e.g. \cite{hl, ho}), such that:
$$
\begin{array}{l}
\E[\integ{0}{T}e^{\alpha
s}\{(\bar{Y}_s-\bar{Y}'_s)^2+|\bar{Z}_s-\bar{Z}'_s|^2+
\integ{E}{}|\bar{V}_s(e)-\bar{V}'_s(e)|^2\lambda(de)\}ds]
\\ \qquad\qquad \qquad \qquad\qquad \leq \bar{C}\E[\integ{0}{T}e^{\alpha
s}\{|y_s -y'_s|^2+|z_s-z'_s|^2+\|v_s-v'_s\|^2\}ds].\end{array}$$
Then the mapping $\Phi$ is a contraction which then has a unique
fixed point $(Y,Z,V)$ which actually belongs to ${\cal S}^{2}\times
{\cal H}^{2,d}\times {\cal L}^2$. Moreover there exists $K^\pm \in
{\cal S}^2$ ($K^\pm_0=0$) such that $(Y,Z,V,K^+,K^-)$ is solution
for the reflected BSDE associated with $(f,\xi,L,U)$. $\Box$
\subsection{Appendix}
Throughout this appendix $(\Omega, {\cal F}, ({\cal F}_t)_{t\leq 1},
P)$ is the same as in Section 2.

Let $\eta:=(\eta_t)_{t\leq 1}$ be a RCLL, ${\cal P}$-measurable
process with values in $R$ and of a class [D]. The Snell envelope of
the process $\eta$, which we denote $R(\eta):=(R(\eta)_t)_{t\leq 1}$
is the lowest RCLL ${\cal F}_t$-supermartingale of class [D] which
dominates $\eta$, i.e., $P-a.s., \,\, \forall t\leq 1,R(\eta)_t\geq
\eta_t$. It has the following expression (see e.g. \cite{el}):
$$ P-a.s.,\,\,
\forall t\leq 1, \,\,R(\eta)_t=\mbox{esssup}_{\tau\geq
t}E[\eta_\tau|{\cal F}_t]\,\,\,(R(\eta)_1=\eta_1).$$

We now give some properties of the Snell envelope of processes.
\medskip

\no \bf{[A1]}: Let $(U^n)_{n\geq 0}$ be a non-decreasing sequence of
${\cal P}$-measurable, RCLL, $R$-valued processes of class [D] which
converges pointwisely to $U$ another RCLL, $R$-valued, ${\cal
P}$-measurable process of class [D], then $P-a.s.$, for any $t\leq
1$, $R(U^n)_t\nearrow R(U)_t$.
\medskip

\noindent $Proof$: Actually for any $n\geq 0$, $P$-a.s. $\forall$
$t\leq 1$, $R(U^n)_t\leq R(U)_t$. Therefore $P-a.s.$, for any $t\leq
1$, $\lim_{n\rightarrow \infty}R(U^n)_t\leq R(U)_t$. Note that the
process $(\lim_{n\rightarrow \infty}R(U^n)_t)_{t\leq 1}$ is an RCLL
supermartingale of class [D] since it is a limit of a non-decreasing
sequence of supermartingales (see e.g. \cite{dm}, pp.86). But
$U^n\leq R(U^n)$ implies that $P-a.s.$, $\forall t\leq 1$, $U_t\leq
\lim_{n\rightarrow \infty}R(U^n)_t$ and then $R(U)_t\leq
\lim_{n\rightarrow \infty}R(U^n)_t$ since the Snell envelope of $U$
is the lowest supermartingale which dominates $U$. It follows that
$P-a.s.$, for any $t\leq 1$, $\lim_{n\rightarrow
\infty}R(U^n)_t=R(U)_t$, whence the desired result. $\Box$
\bigskip

\noindent \bf{[A2]}: \it{Doob-Meyer decomposition of Snell
envelopes}
\medskip

Let $\eta:=(\eta_t)_{t\leq 1}$ be a RCLL, ${\cal P}$-measurable
process with values in $R$ and of a class [D], and
$R(\eta):=(R(\eta)_t)_{t\leq 1}$ its Snell envelope. Then there
exist an RCLL  ${\cal F}_t$-martingale $(M_t)_{t\leq 1}$ and a
non-decreasing RCLL ${\cal F}_t$-predictable process $(K_t)_{t\leq
1}$ ($K_0=0$) such that: $$P-a.s.,\,\, \forall t\leq
1,\,\,R(\eta)_t=M_t-K_t.$$ Moreover we have:

$(i)$ if $R(\eta)$ belongs also to ${\cal S}^2$ then
$E[K_T^2]<\infty$

$(ii)$ if $K^c$ (resp. $K^d$) denotes the continuous (resp. purely
discontinuous) part of $K$ then $K^d$ is ${\cal F}_t$-predictable
and $\{\Delta K^d >0\}\subset \{R(\eta)_-=\eta_-\}$ and
$\Delta_tK^d=(\eta_{t-}-R(\eta)_t)^+1_{\{R(\eta)_{t-}=\eta_{t-}\}}$.
\medskip

\no $Proof$: The existence of $M$ and $K$ is just the Doob-Meyer
decomposition of supermartingales of class [D] (see \cite{dm},
pp.221). Besides if $R(\eta)$ belongs to ${\cal S}^2$ then the
process $K$ is so. This a direct consequence of the dual predictable
projection of $K$ (see \cite{dm}, pp. 221). The proof of \\
$\{\Delta K^d >0\}\subset \{R(\eta)_-=\eta_-\}$ is given in
(\cite{el}, pp.131). Finally since the filtration is generated by a
Brownian motion and an independent Poisson measure the jumps of $M$
occur only at inaccessible stopping times. Therefore when $K^d$
jumps, which is a predictable process, the process $R(\eta)$ has the
same jump. It follows that
$\Delta_tK^d=(R(\eta)_{t-}-R(\eta)_t)^+1_{\{R(\eta)_{t-}=\eta_{t-}\}}=
(\eta_{t-}-R(\eta)_t)^+1_{\{R(\eta)_{t-}=\eta_{t-}\}}$. $\Box$
\medskip

Let $X:=(X_t)_{t\leq 1}$ be a process of class [D]. The predictable
projection of $X$, which we denote by $X^p$, is an ${\cal
F}_t$-predictable process which satisfies $E[X_\tau|{\cal
F}_{\tau-}]=X^p_\tau$ for any predictable stopping time. The process
$X$ is called \it{regular} if it satisfies $X^p_t=X_{t-}$, for any
$t\leq 1$.
\medskip

The following result is related to the existence of an optimal
stopping time when the Snell envelope is $regular$.
\medskip

\no \bf{[A3]} Let $\eta$ be a process of ${\cal S}^2$ and $R(\eta)$
its Snell envelope whose decomposition is $M-K$. For $t\leq 1$, let
$\tau_t$ be the stopping time defined as follows:
$$
\tau_t=\inf\{s\geq t, K_s-K_t>0\}\wedge 1.$$ If $R(\eta)$ is regular
then $K^d\equiv 0$ and $\tau_t$ is optimal after $t$, i.e., it
satisfies:
\medskip

$(i)$ $E[\eta_{\tau_t}]=\sup_{\tau\geq t}E[\eta_\tau]$

$(ii)$ $R(\eta)_{\tau_t}=\eta_{\tau_t}$ and
$(R(\eta)_{s\wedge\tau_t})_{s\geq t}$ is an ${\cal F}_s$-martingale.
\medskip

A word about the proofs of those facts. The continuity of $K$ when
$R(\eta)$ is regular is stated in (\cite{dm}, pp.214). As for the
optimality of $\tau_t$, one can see e.g (\cite{el}, pp. 140). $\Box$
\bigskip

\ni \bf{Acknowledgement}. This paper has been carried out when the
second author visited Universit\'e du Maine (Le Mans, France). Their
hospitality was greatly appreciated. $\Box$ \small{

}

\begin{thebibliography}{99}
\bibitem{bbp} G. Barles, R.\ Buckdahn, E.\ Pardoux: BSDEs and
integral-partial differential equations. \it{ Stochastics 60, pp.
57-83, 1997.}

\bibitem{jeanbl} T.Bielecki, S.Crepey, M.Jeanblanc, M.Rutkowski: Defautable
options in a Markovian intensity model of credit risk.
\it{Mathematical Finance, Vol. 18, Issue 4, pp. 493-518, October
2008 }

\bibitem{CK} J.Cvitanic, I.Karatzas: Backward SDEs with reflection
and Dynkin games, \it{Annals of Probability 24 (4), pp. 2024-2056
(1996)}

\bibitem{dm} C.Dellacherie, P.A.Meyer: Probabilit\'{e}s et Potentiel. \it{Chap.
V-VIII. Hermann, Paris (1980).}

\bibitem{el} N.El-Karoui: Les aspects probabilistes du contr\^{o}le
stochastique, $in$ \it{ Ecole d'\'{e}t\'{e} de Saint-Flour. Lecture
Notes in Mathematics 876, pp. 73-238. Springer Verlag Berlin.}

\bibitem{ekp} N.El-Karoui, C.Kapoudjian, E.Pardoux, S.Peng, M.C.Quenez:
Reflected solutions of backward SDE's and related obstacle problems
for PDE's. \it{Annals of Probability 25 (2) (1997), pp.702-737.}

\bibitem{[EPQ]} N.El-Karoui, E.Pardoux, M.-C.Quenez: Reflected backward
SDEs and American options, \it{ in: Numerical Methods in Finance
(L.Robers and D. Talay eds.), Cambridge U. P., 1997, pp. 215-231}

\bibitem{[H]} S.Hamad\`ene: Mixed Zero-sum differential game and American
game options, \it{SIAM JCO, Vol. 45 (2), pp.496-518 (2006).}

\bibitem{h} S.Hamad\`{e}ne: Reflected BSDE's with discontinuous
barriers and application. \it{ Stochastics and Stochastics Reports
Vol.74 (3-4), pp. 571-596 (2002).}

\bibitem{hl} S.Hamad\`{e}ne, J.P.Lepeltier: Reflected Backward SDE's and
Mixed Game Problems. \it{ Stochastic Processes and their
Applications 85 (2000) pp. 177-188.}

\bibitem{ho} S.Hamad\`{e}ne, Y.Ouknine Reflected Backward stochastic
differential equation with jumps and random obstacle. \it{ EJP Vol.
8 pp. 1-20 (2003).}

\bibitem{iw} N. Ikeda, S. Watanabe: Stochastic Differential Equations and
Diffusion Processes, \it{North Holland/Kodansha (1981).}
















\bibitem{lepeltiermingyu} J.P. Lepeltier, X.Mingyu: Penalization method for
reflected backward stochastic differential equations with one
$r.c.l.l.$ barrier.\it{ Statist. Probab. Lett. 75 (2005), no. 1,
pp.58-66.}

\bibitem{pengmingyu} X.Mingyu, S.Peng: The smallest $g$-supermartingale and
reflected BSDE with single and double L2 obstacles, \it{ Annales de
l'IHP (B), Probability and Statistics, Vol.41, Issue 3, May-June
2005, pp. 605-630.}

\bibitem{pp1} E.Pardoux, S.Peng: Adapted Solutions of Backward Stochastic
Differential Equations. \it{ Systems and Control Letters 14,
pp.51-61, 1990.}

\bibitem{[Pa]}  E. Pardoux: BSDEs, weak convergence and homogenization of
semilinear PDEs, \it{in: F. Clarke and R. Stern (eds), Nonlin.
Anal., Dif. Equa. and Control, pp. 503-549 (1999), Kluwer Acad.
Publi., Netherlands}

\bibitem{[PP1]} E. Pardoux, S. Peng: Backward Stochastic Differential
Equations and Quasilinear Parabolic

Partial Differential Equations, \it{in: Stochastic Differential
Equations and their Applications (B. Rozovskii and R. Sowers, eds.),
Lect. Not. Cont. Inf. Sci., vol.176, Springer, 1992, pp. 200-217}

\bibitem{monlim} Peng, S. (1999): Monotonic limit theory of BSDE and nonlinear decomposition
theorem of Doob-Meyer's type, \it{ Probability Theory and Related
Fields, 113, pp. 473-499.}



\bibitem{tl} S. Tang and X. Li: Necessary condition for optimal control of
stochastic systems with random jumps, \it{SIAM JCO 33 2, pp.
1447-1475, (1994).}
\end{thebibliography}
\end{document}